\documentclass[12pt]{article}
\usepackage{amsgen,amsmath,amstext,amsbsy,amsopn,amsfonts,amssymb}
\newtheorem{theorem}{Theorem}
\newtheorem{lemma}[theorem]{Lemma}
\newtheorem{corollary}[theorem]{Corollary}
\newtheorem{proposition}[theorem]{Proposition}
 \newtheorem{defi}[theorem]{Definition}
\newenvironment{definition}{\begin{defi}\rm}{\end{defi}}
\newtheorem{exa}[theorem]{Example}
\newenvironment{example}{\begin{exa}\rm}{\end{exa}}
\newtheorem{rem}[theorem]{Remark}
\newenvironment{remark}{\begin{rem}\rm}{\end{rem}}
\newtheorem{rems}[theorem]{Remarks}

\newtheorem{ack}[theorem]{Acknowlegment}

\def\bsq{\blacksquare\medskip}

\def\H{\mathcal H}
\def\L{\mathcal L}
\def\K{\mathcal K}

\def\Hpi{(\pi, {\mathcal H})}

\def\ZZ{{\mathbb Z}}
\def\CCC{{\mathbb C}}
\def\RRR{{\mathbb R}}

\def\RR+{{\mathbb R}^*}

\def\HS{\rm HS}
\def\Q_p{{\mathbb Q}_p}
\def\S1{{\mathbb S}^1}

\def\PT{Property (T)}
\def\PTB{Property (${\rm T}_B$)}
\def\eps{\varepsilon}
\def\Ga{\Gamma}
\def\ga{\gamma}
\def\vfi{\varphi}

\begin{document}

\title{A  \PT\ for $C^*$-algebras  }
\author{ Bachir Bekka }
\date{\today }

\maketitle

\begin{abstract}
We define a notion of Property (T) for an arbitrary
 $C^*$-algebra $A$ admitting a tracial state.
We extend this to a notion
of Property (T)  for the pair $(A,B),$ 
where $B$ is a  $C^*$-subalgebra of $A.$ 
Let $\Ga$ be a discrete
group and $C^*_r(\Ga)$ its reduced algebra. 
We show that 
$C^*_r(\Ga)$ has \PT\ if and only if
the group $\Ga$ has \PT.
More generally, given a subgroup
$\Lambda$ of $\Gamma$, the pair
$(C^*_r(\Ga),C^*_r(\Lambda)) $ 
has \PT\ if and only if the
pair of groups  $(\Ga, \Lambda)$ has \PT.
\end{abstract}

\section{Introduction}
\label{S1}
Let $G$ be a locally compact group. Recall that
$G$ has Kazhdan's \PT\ if, whenever a unitary representation
$\Hpi$ of $G$ almost
has invariant vectors, $\H$ has a non-zero invariant vector
(see  \cite{Kazhdan}; \cite{HaVa}).

A. Connes defined a notion of Property (T) for type $II_1$-factors
\cite{Connes} and it was
proved  in  \cite{CoJo} that a discrete ICC-group
$\Ga$ has \PT\ if and only if 
 the von Neumann algebra  $L(\Ga)$ generated
by its left regular representation has \PT.
Jolissaint \cite{Jolissaint} gave an extension of
this  result  to arbitrary discrete groups.
(Recall that a group is ICC if all its conjugacy
classes, except the unit element,
are infinite.)
For an extensive study of \PT\  
for von Neumann algebras and, more generally, for inclusions
of von Neumann algebras,
see the notes of S. Popa \cite{Popa}. 

We introduce \PT\ for
an arbitrary unital $C^*$-algebra 
admitting a tracial state  (see
Remark~\ref{Rem1} below) through an adaptation
of
Connes' definition. This allows 
us to express
\PT\ for an arbitrary  countable group in terms of its {reduced}
$C^*$-algebra.

Let $A$ be a $C^*$-algebra or a von Neumann algebra.
 A Hilbert {bimodule}
over $A$ is a Hilbert space $\H$ carrying two commuting 
representations, one of  $A$ and one of the opposite algebra $A^0;$ we will 
write 
$$\xi\mapsto x\xi y,  \qquad \forall \xi\in \H,\quad \forall x,y\in A.$$ 
(In the case where $A$ is a von Neumann algebra,
we consider only {normal} representations of $A.$)

Recall that a tracial state
on a unital $C^*$-algebra $A$
is a positive linear functional
$\tau: A\to \CCC$ such that 
$\tau (xy)=\tau(yx)$ for all
$x,y\in A$ and $\tau(1)=1.$   
\begin{definition}
\label{Def1} 
Let $A$ be either a unital $C^*$-algebra admitting a tracial state
or a finite von Neumann  algebra.
We say that $A$ has \PT\ if there
exist a finite subset $F$ of $A$ and $\eps>0$
such that the following property  holds:
if a Hilbert bimodule $\H$ for  $A$ contains a unit vector 
$\xi$ which is $(F,\eps)$-central,
that is, such that
 $$ \Vert y\xi-\xi y\Vert<\eps, \quad \forall y\in F,$$
then $\H$ has a non-zero central vector,
that is, a vector $\eta\in \H$
such that $x\eta =\eta x$ for all $x\in A.$
\end{definition}

Our main result  shows
that \PT\ for a countable discrete group $\Ga$
only depends  on its reduced $C^*$-algebra
$C^*_r(\Ga).$ Recall that
$C^*_r(\Ga)$ is the norm closure
of the linear span of 
$\{\lambda_{\Ga} (\ga)\ :\ \ga\in \Ga\}$
in $\L(\ell^2(\Ga)),$
where $\lambda_{\Ga}$ is the left regular representation
of $\Ga.$ The mapping
$x\mapsto \langle x\delta_{e},\delta_e\rangle$
is a (faithful) tracial  state on $C^*_r(\Ga).$
The von Neumann algebra $L(\Ga)$ of  $\Ga$ is the
von Neumann subalgebra of $\L(\ell^2(\Ga))$
generated by $C^*_r(\Ga),$ in other words, 
$L(\Ga)$ is the closure of the linear span of 
$\{\lambda_{\Ga} (\ga)\ :\ \ga\in \Ga\}$ in the weak 
operator topology.
We shall also consider the maximal (or full) $C^*$-algebra
$C_{\max}^*(\Ga)$ of $\Ga,$ that is,
the completion of the group algebra
$\CCC \Ga$ with respect to the
norm
$$
\Vert \sum_{\ga\in \Ga}{a_{\ga} \ga}\Vert =\sup_{\pi} \Vert\sum_{\ga\in \Ga}{a_{\ga} \pi(\ga)}\Vert,
$$
where $\pi$ runs over the set of all equivalence classes of 
cyclic unitary representations
of $\Ga.$   

\begin{theorem}
\label{Theo1}
Let $\Ga$ be a countable discrete group, and let
$A$ be a $C^*$-algebra. Assume that $A$
  is a quotient of 
$C_{\rm max}^*(\Ga)$ and that $C^*_r(\Ga)$ is a quotient
of $A.$ 
The following  properties are equivalent:
\begin{itemize}
\item[(i)] $\Ga$ has \PT;
\item[(ii)] $A$ has \PT;
\item[(iii)] $L(\Ga)$ has \PT;
\end{itemize}
\end{theorem}
Observe that a $C^*$-algebra $A$ as in theorem above
has at least
one tracial state. Indeed, since  $C_{r}^*(\Ga)$ is
a quotient of $A,$ the canonical trace
of $C_{r}^*(\Ga)$ lifts to a tracial state
of $A.$ 

Apart from
$C_{r}^*(\Ga)$ and $C_{\rm max}^*(\Ga),$  examples of 
$C^*$-algebras as in the theorem
arise from  unitary representations of $\Ga$
which vanish at infinity.
If $\Hpi$ is a unitary representation
of the discrete group $\Ga,$ we denote by 
$C^*_\pi (\Ga)$  the
$C^*$-subalgebra of $\L(\H)$ generated
by $\{\pi(\ga)\ :\ \ga\in \Ga\}.$
Such an algebra is of course a quotient of 
$C_{\rm max}^*(\Ga).$ 
The following  is Proposition~1
in \cite{BeHa} (see also the comments
at this end of this paper).
\begin{proposition}
\label{PropC0}
{\textbf (\cite{BeHa})}
 Let $\Ga$ be an ICC group, and let $\pi$ be a unitary 
representation of $\Ga.$ Assume that
all matrix coefficients of $\pi$ are in $c_0(\Ga).$
Then $C^*_r (\Ga)$ is a quotient of
$C^*_\pi(\Ga).$ 
\end{proposition}
Using \cite[Th\'eor\`eme~1]{BeHa}
and using the fact that a lattice $\Ga$ in a locally compact 
group $G$ has \PT\ if and only if  $G$ has \PT,
we obtain
the following class of examples.
\begin{corollary}
\label{CorC0}
Let $G$ be a connected, non compact simple real group
with trivial centre.
Let $\pi$ be a  
unitary representation of $G$ 
which is not a multiple of the trivial representation.
For any lattice  $\Ga$ 
in $G,$ the following properties are equivalent:
\begin{itemize}
\item[(i)] $G$  has \PT;
\item[(ii)] $C^*_{\pi\vert_\Ga} (\Ga)$  has \PT.
\end{itemize} 
\end{corollary}
Observe that, if in the previous
corollary the representation
$\pi$ is irreducible and not square-integrable, then
$C^*_{\pi\vert_\Ga}(\Ga)$ is a primitive $C^*$-algebra.
Indeed, as shown by  Cowling and Steger
\cite{CoSt}, in this case, the restriction
of $\pi$ to $\Ga$ is irreducible.
The simple real Lie groups with \PT\ are 
known: these are exactly the simple real Lie 
groups which are not locally isomorphic to $SO(n,1)$
or $SU(n,1)$ (see \cite{HaVa}). 

We now define \PT\ for pairs consisting of a $C^*$-algebra
and a $C^*$-subalgebra. Recall that the notion of \PT\ for groups can be generalized
to a notion of \PT\ for pairs $(G,H),$ also called relative \PT, 
where $G$ is a locally compact group  and $H$ a closed
subgroup. Such a pair  $(G,H)$
 is said to have \PT\ if, whenever a unitary representation
$\Hpi$ of $G$ almost
has invariant vectors, it has a non-zero $H$-invariant vector
(see \cite{HaVa}). It is clear that \PT\ for the group
$G$ is equivalent to
\PT\ for the pair $(G,G).$ A prominent example
of a pair with \PT\ is the pair $(SL_2(\ZZ)\ltimes\ZZ^2, \ZZ^2),$
where $SL_2(\ZZ)\ltimes\ZZ^2$ is the semi-direct product
for the natural action of $SL_2(\ZZ) $ on $\ZZ^2.$

\begin{definition}
\label{Def2} 
Let $A$ be  a unital $C^*$-algebra admitting a tracial state
(respectively, a finite von Neumann  algebra) and let $B$ be  
a $C^*$-subalgebra  (respectively, a von Neumann  subalgebra)
of $A.$
The pair $(A,B)$ has \PT\ if there
exist a finite subset $F$ of $A$ and $\eps>0$
with the following property:
if a Hilbert bimodule $\H$ for  $A$ contains a unit vector 
$\xi$ which is $(F,\eps)$-central,
then $\H$ has a non-zero $B$-central vector,
that is, a vector $\eta\in \H$
such that $b\eta =\eta b$ for all $b\in B.$
\end{definition}
Theorem~\ref{Theo1} can be extended to pairs as follows.
Let $\Ga$ be a discrete group  and $\Lambda$ a subgroup
of $\Ga.$  The restriction of
the regular representation $\lambda_\Ga$ of
$\Ga$  to $\Lambda$ is a multiple of $\lambda_\Lambda.$
This implies that the canonical embedding of group algebras
$\CCC\Lambda\to \CCC\Ga$ extends to  isometric embeddings
 $C_{r}^*(\Lambda)\to C_{r}^*(\Ga)$ and $L(\Lambda)\to L(\Ga).$ 
Moreover, since every function
of positive type on $\Lambda$ extends (in a trivial way)
to  a function of positive type on $\Ga,$ every unitary representation of $\Lambda$ is 
contained in the restriction of a unitary representation of $\Ga.$
This shows that the embedding $\CCC\Lambda\to \CCC\Ga$
extends also to an isometric embedding $C_{\max}^*(\Lambda)\to C_{\max}^*(\Ga).$

\begin{theorem}
\label{Theo2}
Let $\Ga$ be a countable discrete group
and $\Lambda$ a subgroup of $\Ga.$
The following  properties are equivalent:
\begin{itemize}
\item[(i)] the pair $(\Ga, \Lambda)$ has \PT;
\item[(ii)] the pair $(C_{\max}^*(\Ga),C_{\max}^*(\Lambda) )$ has \PT;
\item[(iii)] the pair $(C_{r}^*(\Ga),C_{r}^*(\Lambda) )$ has \PT;
\item[(iv)] the pair $(L(\Ga), L(\Lambda))$ has \PT.
\end{itemize}
\end{theorem}

\begin{remark}
\label{Rem 5}
The definition given above
of \PT\ for a pair
of von Neumann algebras 
appeared in the first version of \cite{Popa03}.
In  the final version of  \cite{Popa03},
this definition was replaced by 
a stronger notion of \PT,
called ``rigid embedding" of von Neumann algebras.
This  notion  seems to be the more useful one
(compare also \cite {PePo}). Moreover,  a proof
of the equivalence of (i) and \PT\ for the
pair  $(L(\Ga), L(\Lambda))$ in the 
strong  sense was given in \cite[5.1 Proposition]{Popa03}.
\end{remark}

For the proof of Theorem~\ref{Theo2}, we will use the following
characterization of \PT\ for pairs of groups,
which is of independent interest.
This result is an extension of \cite[Theorem 1]{BeVa},
where the case $H=G$ is considered.

\begin{theorem}
\label{Theo-BeVa}
 Let $G$ be a $\sigma$-compact locally compact group
and $H$ a closed subgroup of $G.$
The following properties
are equivalent:
\begin{itemize}
\item[(i)]  The pair $(G,H)$ has \PT;
\item[(ii)] if a unitary representation $\Hpi$
of $G$ almost has invariant vectors, then 
$\H$   contains  a non-zero finite dimensional 
subspace which is invariant under $H.$
\end{itemize}
\end{theorem}

This article is organized as follows.
In Section~\ref{S:BeVa}, we give the proof
of Theorem~\ref{Theo-BeVa}.
The proof given in \cite{BeVa}  carries over to pairs,
but, for the convenience of the reader, we prefer to give
a different and shorter  proof of this result
(for another proof, see also \cite[Lemma 2.1]{Jolissaint2}).
Section~\ref{S:TheoPrincipal} 
is devoted to the proof
of Theorem~\ref{Theo2}.
In Section~\ref{S:Nuclearity},
the relationship between \PT\ and nuclearity
is explored. We give there a $C^*$-algebraic analogue
(see Proposition~\ref{Proposition} below)
of the known fact that a discrete group which
has \PT\ and which is amenable has to be finite.
Section~\ref{S:Rem} contains some remarks
around our definition of \PT.

\noindent 
\textbf{ Acknowledgments.}
It is a pleasure to thank E. Blanchard, 
M. Gromov,  P. de la Harpe,
G. Skandalis and T. Steger for helpful comments.

\section{Proof of Theorem \ref{Theo-BeVa}}
\label{S:BeVa}

 It is obvious that (i) implies (ii). To show the converse,
assume that $(G,H)$ does not have \PT. 
Since $G$ is $\sigma$-compact,  there exists  a
function $\psi:G\to\RRR$ of conditionally negative type
such that the restriction $\psi\vert _H$ is unbounded.
This is the Delorme-Guichardet theorem
(see\cite[Theorem 5.20]{HaVa}) which, as
is easily seen, carries over to pairs.
We can clearly assume that $\psi(e)=0.$ 
By Schoenberg's theorem, $\vfi_t=\exp(-t\psi)$ is a normalized function of positive
type on $G$ for every $t>0.$
Let $(\pi_t, \H_t, \xi_t)$ be the triple
associated to $\vfi_t$ by GNS-contruction, that is,
$\pi_t$ is a cyclic unitary representation of
$G$ on the Hilbert space $\H_t,$
with cyclic unit vector $\xi_t$ such that
$\vfi_t(g)=\langle\pi_t(g)\xi_t,\xi_t\rangle$ for all
$g\in G.$

Let $(t_n)_n$ be a sequence of positive numbers with
$\lim_n t_n=0$ and set
 $$\pi=\bigoplus_n \pi_{t_n}.$$
Since $\lim_{n} \exp(-t_n\psi(g))=1$ uniformly on
compact subsets of $G,$
the representation $\pi$ almost has invariant vectors.
We claim that, for every $t>0,$
the restriction of the representation $\pi_t$ to $H$
and hence  the restriction of  $\pi$
of $H$  contains  no (non-zero) finite 
dimensional subrepresentation. 

Fix $t>0.$  The GNS-construction
applied to $\psi_t=t\psi$ shows that there exists an action $\alpha_t$
of $G$ by affine isometries on some real  Hilbert space 
such that $\psi_t(x) =\Vert \alpha_t(x)0\Vert^2$ for all $x\in G.$
Since $\psi_t\vert_H$ is unbounded, we can find a sequence
$(x_i)_i$  in $H$ such that $
\lim_i\psi_t(x_i)=+\infty.$
For $a,b\in G,$ we have 
\begin{eqnarray*}
\sqrt{\psi_t(ab)}&=&\Vert \alpha_t(ab)0-0\Vert
=\Vert \alpha_t(b)0-\alpha_t(a^{-1})0\Vert\\
&\leq& \Vert \alpha_t(b)0-0\Vert +\Vert 0- \alpha_t(a^{-1})0\Vert
=\Vert \alpha_t(b)0-0\Vert +\Vert  \alpha_t(a)0-0\Vert\\
&=&\sqrt{\psi_t(a)} +\sqrt{\psi_t(b)}.
\end{eqnarray*}
Hence, we have
\begin{eqnarray*}
\sqrt{\psi_t(ax_ib)}
&\leq&\sqrt{\psi_t(x_i)}+\sqrt{\psi_t(b)}+\sqrt{\psi_t(a)}.
\end{eqnarray*}
and therefore $\lim_i\psi_t(ax_ib)=+\infty,$ 
that is, 
$$
\lim_i\varphi_t(ax_ib)=0 , \qquad  \forall a,b\in G.
$$
Since $\xi_t$ is a cyclic vector for $\pi_t,$  it follows that
$\lim_i\langle \pi_t(x_i)\eta_1,\eta_2\rangle =0$
for all $\eta_1, \eta_2\in \H_t.$
Let $\K$ be a finite
dimensional $H$-invariant 
subspace of $\H_t$ and let $\eta\in \K.$
Upon passing to a subsequence,
we can assume that $(\pi_t(x_i)\eta)_i$  is \emph{norm} convergent
to some vector $\eta_0\in \K.$ 
Since $\lim_i\langle \pi_t(x_i)\eta,\eta_0 \rangle=0,$
we have $\eta_0=0$ and hence $\eta=0.$ 
Therefore $\K=0.$
$\bsq$

\section{Proof of Theorem~\ref{Theo2}} 
\label{S:TheoPrincipal}
 Let $G\to H$ be a continuous homomorphism
between locally compact groups
with dense image. It is well-known
that if $G$ has \PT, then
 $H$ has \PT (see \cite{HaVa}).
The corresponding statement for $C^*$-algebras is as follows
and its proof is straightforward.
\begin{lemma}
\label{Lem1}
 Let $A_1$  be a unital 
$C^*$-algebra admitting tracial states
and $B_1\subset A_1$ be a  $C^*$-subalgebra.
Let $A_2$ be  a unital 
$C^*$-algebra admitting tracial states
(respectively, a finite von Neumann algebra), and 
let $f: A_1\to A_2$ be a unital $*$-homomorphism.
Let $ B_2=f(B_1)$ (respectively, let $B_2$
be the von Neumann subalgebra of $A_2$
generated by $f(B_1)$).
  If the
pair $(A_1, B_1)$ has \PT, then the pair $(A_2, B_2)$ has \PT.
$\bsq$
\end{lemma}
We proceed now with the proof of 
Theorem~\ref{Theo2}.
To show that (i) implies (ii), assume that the pair
$(\Ga, \Lambda)$ has \PT. Then there exists
a Kazhdan pair $(S,\eps),$ that is,
a finite subset $S$ of $\Ga$ and $\eps>0$
with the following property:  
if $\Hpi$ is a unitary representation
of $\Ga$ and if there exists a unit vector
$\xi\in \H$ with $\max_{s\in S} \Vert\pi(s)\xi-\xi\Vert<\eps,$
then $\H$ has a non-zero $\Lambda$-invariant vector.

Let $\H$ be a Hilbert bimodule of $C_{\max}^*(\Ga).$
Viewing $\Ga$
as a subset of $C_{\max}^*(\Ga),$
we define two commuting unitary representations $\pi_1$ and $\pi_2$
of $\Ga$ on $\H$  by 
$$
\pi_1(\ga)\xi = \ga \xi\qquad\text{and}\qquad \pi_2(\ga) \xi =\xi \ga^{-1}
$$
for all $\ga\in \Ga$ and $\xi\in \H.$
Assume that $\H$ has a unit vector $\xi$
such that 
$$ \Vert s\xi-\xi s\Vert<\eps, \quad \forall s\in S.$$
Then 
$$ \Vert \pi_1(s)\pi_2(s)\xi-\xi \Vert<\eps, \quad \forall s\in S.$$
Hence, there exists a non-zero 
vector $\eta\in \H$
with $$\pi_1(\ga)\pi_2(\ga)\eta=\eta, \qquad \forall\ga\in \Lambda.$$
Since the linear span of $\Lambda$ is dense in $C_{\max}^*(\Lambda),$
we have $x\eta=\eta x$ for all $x\in C_{\max}^*(\Lambda),$
that is, $\eta$ is $C_{\max}^*(\Lambda)$-central.

The fact that  (ii) implies (iii) and that
(iii) implies (iv) follows from the previous lemma,
since $C_{r}^*(\Ga)$ is a quotient of $C_{\max}^*(\Ga)$
and since $C_{r}^*(\Ga)$ is weak-* dense in $L(\Ga).$

The proof  that (iv) implies (i)
is an adaptation of
the proof of \cite[Theorem~2]{CoJo}. Indeed,
assume that $(L(\Ga),  L(\Lambda))$
has \PT. Choose a finite subset $F$
of  $C_{r}^*(\Ga)$ and $\eps>0$ as in
Definition~\ref{Def2}.
We may assume that $\Vert y\Vert \leq 1$ for all
$y\in F.$ Let $S$  be a finite
subset of $\Ga$ such that
$$\sum_{\ga \in \Ga\setminus S}|(y\delta_{e})(\ga)|^2 <\eps^2/9$$
for all $y\in F.$

Let $\pi$ be a
unitary representation of $\Ga$.
Define two commuting unitary representations $\pi_1$ and $\pi_2$
of $\Ga$ on the Hilbert space tensor product
$\ell^2(\Ga)\otimes \H=\ell^2(\Ga, \H)$ by 
$$
\pi_1(\ga)\widetilde\xi(x) =  \widetilde\xi(\ga^{-1}x)
$$
and 
$$
\pi_2(\ga)\widetilde\xi(x) =  \pi(\ga)\widetilde\xi(x \ga)
$$
for all $\ga\in \Ga$ and $\widetilde\xi\in \ell^2(\Ga, \H)$.
Since $\pi_1$ and $\pi_2\simeq \lambda_{\Ga}\otimes \pi$
are equivalent to multiples of the regular representation
$\lambda_{\Ga},$ they extend to commuting   representations 
of $L(\Ga),$ so that $\ell^2(\Ga, \H)$
is a bimodule of $L(\Ga).$

Assume now that there exists a unit vector
$\xi\in \H$ such that  
$$\max_{s\in S} \Vert\pi(s)\xi-\xi\Vert<\eps/3.$$
Let $\widetilde\xi\in \ell^2(\Ga, \H)$
 be defined by $\widetilde\xi (e)= \xi$
and $\widetilde\xi (x)= 0$ otherwise,
that is, $\widetilde\xi=\delta_{e}\otimes \xi.$
Then, for every $y\in F,$
we have
\begin{eqnarray*}
\Vert y\widetilde\xi-\widetilde\xi y\Vert^2&=&
\sum_{\ga\in \Ga} |y\delta_{e}(\ga)|^2\Vert \pi(\ga)\xi-\xi\Vert^2\\
&\leq& \frac{4}{9}\eps^2 +\sum_{\ga\in S} |y\delta_{e}(\ga)|^2\Vert \pi(\ga)\xi-\xi\Vert^2\\
&\leq& \frac{4}{9}\eps^2 +\frac{1}{9}\eps^2 <\eps^2.
\end{eqnarray*}

Hence, there exists a non-zero  vector $\eta$ in $\ell^2(\Ga, \H)$
which is $L(\Lambda)$-central.
We then have
$$\eta (\ga^{-1} x\ga)=\pi(\ga)\eta (x), \qquad \forall x\in \Ga, \ \ga\in \Lambda.$$
 In particular, $x\to \Vert \eta(x)\Vert$
is a  non-zero function in $\ell^2(\Ga)$ which is invariant
under conjugation by elements from $\Lambda.$
Let $x_0\in \Ga$ be such that
$\eta (x_0)\neq 0.$  It follows that 
its $\Lambda$-conjugacy class
$\{\ga^{-1} x_0\ga\ :\ \ga\in \Lambda\}$ is finite.
Let $\xi_0=\eta(x_0)\in \H.$
Then $\{\pi(\ga)\xi_0 \ :\ \ga\in \Lambda\}$
is a finite subset of $\H$ and its linear span defines
a finite dimensional subrepresentation of $\pi\vert_\Lambda.$
It follows from
Theorem~\ref{Theo-BeVa} that the pair $(\Ga, \Lambda)$ has \PT.
$\bsq$

\begin{remark}
\label{RemTheo1}
Let $A$ be a $C^*$-algebra which 
  is a quotient of 
$C_{\rm max}^*(\Ga)$ and such that $C^*_r(\Ga)$ is a quotient
of $A.$ The same proof as above (with $\Lambda=\Ga$
and with $A$ instead of $C^*_r (\Ga)$)
shows that $\Ga$ has \PT\ if and only
if $A$ has \PT. This proves Theorem~\ref{Theo1}.
\end{remark}

\section{\PT\ and nuclearity}
\label{S:Nuclearity}

As is well-known, a locally compact group with \PT\
which is  amenable 
has to be compact (see \cite[Chap.1, Proposition~6]{HaVa}). 
We will show that a similar fact is true
for $C^*$-algebras with \PT\ 
which are nuclear.
Recall that a $C^*$-algebra $A$ is nuclear
if, for any other $C^*$-algebra $B,$
there a unique pre-$C^*$-norm
on the algebraic tensor product
$A\odot B$ (see \cite[Chap. XV]{Takesaki}).
 Recall also
that, by work
of Connes and Haagerup,
the class of nuclear 
$C^*$-algebras coincides
with the class of $C^*$-algebras
which are amenable in the sense
of B. Johnson (see \cite{Runde}).

Let $\tau$ be  a tracial state 
on the unital $C^*$-algebra $A.$ By the GNS\--con\-struction,  
$\tau$ defines an $A$- Hilbert bimodule, denoted
by $L^2(A,\tau),$ which has a unit central 
vector $\eta\in L^2(A,\tau) $  such that
$\tau(x)=\langle x\eta,\eta\rangle$
for all $x\in A.$ More precisely,
$N=\{x\in A\ :\ \tau(x^*x)=0\}$ is a two-sided
$*$-ideal in $A.$ Define an inner product
on the quotient $A/N$ by 
$\langle x+N, y+N\rangle= \tau(y^*x).$
Let $L^2(A,\tau)$ be the Hilbert space
completion of $A/N.$
For each $a\in A,$ the mappings $x\mapsto ax+N$
and $x\mapsto xa+N$ extend to bounded operators
$\lambda(a)$ and $\rho(a)$ on $L^2(A,\tau),$
defining an $A$-bimodule structure on $L^2(A,\tau).$ 
The von Neumann algebras $\lambda (A)''$ and $\rho(A)''$
are commutant to each other.
Moreover,  $\eta=1+N$ is a cyclic
and separating vector for $\lambda (A)''$ and $\rho(A)''$
and   we have
$\tau(x)=\langle x\eta,\eta\rangle$
for all $x\in A.$

The following proposition is inspired by 
\cite[Proposition 1.2.4.ii]{Popa}.
\begin{proposition}
\label{Proposition}
Let $A$ be a  unital
$C^*$-algebra with \PT. Assume that
$A$ is nuclear.
Then, for any tracial state 
$\tau$ on $A,$  the representation
$\lambda$ on the Hilbert space
$L^2(A,\tau)$ is completely atomic, that is,
$L^2(A,\tau)$ decomposes  as a direct sum of
finite dimensional $A$-submodules.
\end{proposition}

\begin{proof}
Let $L^2(A,\tau)_{\rm fin}$ be the purely atomic part
of $L^2(A,\tau),$ that is, $L^2(A,\tau)_{\rm fin}$ is 
the closed subspace generated by all
finite dimensional $\lambda(A)$-invariant
subspaces.
Let $\H$ be the orthogonal complement
of $L^2(A,\tau)_{\rm fin}.$
Assume, by contradiction, that 
$\H\neq 0.$ We have
$\H= p L^2(A,\tau)$ for some  projection
$p$ belonging to the centre $\lambda (A)''\cap \rho(A)''$
of the finite von Neumann algebra $M=\lambda (A)''.$
The space  $\HS (\H)\simeq\H\otimes \H$ of Hilbert-Schmidt
operators on $\H$ is a $pM$-bimodule
for the action
$$
T\mapsto xTy, \qquad T\in {\HS} ({\H}),\ x, y\in pM.
$$
As $A$ is nuclear,
the finite von Neumann algebra $pM$ acting on $\H$ is injective.
Hence, there exists a conditional expectation
$E: \L(\H)\to pM.$ The state
$\Phi$ on $\L(\H)$
defined by $ \Phi(T)= \tau(E(T))$
is a hypertrace: we have
$$
\Phi(xT)=\Phi(Tx), \qquad \forall T\in\L(\H),\, x\in pM.
$$
It follows  that  there exist almost  central vectors
in $\HS (\H),$ that is,
there exists a net $(T_i)\in \HS (\H)$ with
$\Vert T_i\Vert_{\rm HS}=1$ such that
$$
\lim_i \Vert xT_i-T_ix\Vert_{\rm HS}=0, \qquad \forall x\in pM
$$
(see, for instance,  \cite[Chap. XV, Lemma 3.9]{Takesaki}).
Since $A$ has \PT, we find a non-zero operator $T\in \HS (\H)$
with $xT=Tx$ for all $x\in pM.$
Looking at the eigenspaces of the compact self-adjoint
operator $T^*T,$ we deduce that $\H$ contains a
non-zero finite dimensional $\lambda(A)$-invariant
subspace. This is a contradiction.
$\bsq$
\end{proof}
\begin{corollary}
\label{Cor} 
Let $A$ be a  unital
$C^*$-algebra with \PT. Assume that
there exists a tracial state
$\tau$ of $A$ such that
$L^2(A,\tau)$ is not completely atomic.
Then $A$ is not nuclear.$\bsq$
\end{corollary}

\begin{example}
\label{Exa} 
(i) Let $\Ga$ be a discrete group.
For the canonical trace $\tau$ on $M=C^*_r(\Ga),$
the Hilbert space $L^2(M,\tau)$ can be identified
with $ \ell^2(\Ga),$ with bimodule structure
coming from the left and right representations.
If $\Ga$ is infinite, then $ \ell^2(\Ga)$ has no non-zero finite dimensional
left invariant subspace.
On the other hand,  $C^*_r(\Ga)$ is nuclear
if and only if $\Ga$ is amenable (see
\cite[Theorem 1.1]{Lance}).
Thus, using Theorem~\ref{Theo1}, we recover the fact 
that an infinite discrete  group with \PT\ is not 
amenable.

(ii) All abelian unital $C^*$-algebras are of 
the form $C(X)$ for a compact topological space
$X.$ Such algebras are nuclear.
If $X$ is uncountable,  there exists a
regular probability measure $\mu$ on the
Borel subsets of $X$ which has a
non atomic part.
For such a measure,  $L^2(X,\mu)$ does not decompose
as  a direct sum of finite dimensional
subspaces, under the action
of $C(X)$ as multiplication operators.
As a consequence, $C(X)$ does not have \PT\
if $X$ is uncountable.
 \end{example}

\begin{remark}
\label{Rem9} 
Let $\Ga$ be a discrete group
with \PT\ and $M$ a finite injective  von Neumann
algebra. In connection with Proposition~\ref{Proposition},
the following result of
 Robertson (\cite{Guyan}) is worth mentioning: 
if $\pi: \Ga\to U(M)$ is  a homomorphism
of $\Ga$ into the unitary group of $U(M),$
then  $\pi(\Ga)$ has a compact closure in $U(M)$
for the strong topology; see also \cite{Alain}.
In fact, the result is true more generally
if $M$ has the Haagerup approximation property
(\cite{Guyan}). This is the case for instance
if $M=L(F_2),$ where $F_2$ is the free group
on two generators.
\end{remark}

\section{Some remarks}
\label{S:Rem}

\begin{remark}
\label{Rem1}
Our definition of \PT\ (see Definition~\ref{Def1} above)
 makes sense
for any $C^*$-algebra. However,
with this definition, every unital
$C^*$-algebra $A$ without tracial states
has \PT. Indeed, otherwise,
for every finite subset $F$ of $A$
and every $\eps >0,$
we find a Hilbert bimodule $\H_{F,\eps}$ for  $A$
which has a $(F,\eps)$-central unit  
vector $\xi_{F,\eps}.$ 
Let $\H$ be the direct sum of 
the $\H_{F,\eps}$'s.
Then $\H$ is an $A$-bimodule in an obvious
way.
For every pair $(F,\eps),$
consider the vector state $\varphi_{F,\eps}$
on the algebra $\L (\H)$ 
of all bounded operators on $\H$
defined by $\xi_{F,\eps}\in \H.$ 
Let $\varphi$ be a weak*-limit point of the
net $(\varphi_{F,\eps})_{F,\eps}.$
Then $\tau:A \to \CCC$ defined 
by  $\tau(x)=\varphi (\pi(x))$
is a tracial state on $A,$
where $\pi$ is the representation
on $\H$ given by, say, the left action
of $A.$ This is a contradiction.
\end{remark}

\begin{remark}
\label{Rem 2}
 The definition of \PT\ given
in Definition \ref{Def1} 
makes  sense for an arbitrary 
(normed) algebra $A.$ 
It can  even be further extended
as follows. 
Let $B$ be a fixed Banach
space and consider   all 
(continuous) $A$-bimodule
structures on $B.$ We 
say that $A$ has  \PTB\ if
there exists a finite subset $F$ of $A$ and $\eps>0$
such that, whenever, 
for some $A$-bimodule structure on $B,$ 
there exists a unit $(F,\eps)$-central
 vector,  then $B$ has a non-zero central vector.
\end{remark}

\begin{remark}
\label{Rem 3} A notion of  co-rigid
inclusion $B\subset M$ for a pair 
consisting of a  finite von Neumann algebra
$M$ and a von Neumann subalgebra $B$
was defined 
by Popa (see \cite[Definition 4.1.3]{Popa}, where the inclusion
$B\subset M$ is called rigid; see also \cite{Claire}). 
In a similar way, let $B\subset A$ be a pair 
consisting of a  unital  $C^*$-algebra $A$ admitting a tracial state
 and a $C^*$-subalgebra $B.$ We say that 
$B\subset A$ is co-rigid if there
exist a finite subset $F$ of $A$ and $\eps>0$
such that the following property  holds:
if a Hilbert bimodule $\H$ for  $A$ contains a unit vector 
$\xi$ which is $(F,\eps)$-central and 
which is central for $B$
(that is,  $b\xi =\xi b$ for all $b\in B$)
then $\H$ has a non-zero  vector
which is central for $A.$
\end{remark}

\begin{remark}
\label{Rem 7}
In Definition~\ref{Def1}, we do not require
that the $C^*$-algebra $A$ has a \emph{faithful}
tracial state, that is, a  tracial state $\tau$ such that
$\tau(x^*x)\neq 0$ for all $x\in A$ with $x\neq 0.$
Indeed, 
the group $SL_3(\ZZ)$ has \PT\
but, as was shown by the author (unpublished),
 its maximal $C^*$-algebra $C^*_{\rm max}(SL_3(\ZZ))$
has no faithful  tracial state.
\end{remark}
\begin{remark}
\label{rem8}
Let $M$ be a type $II_1$ factor with \PT.
Let
${\rm Aut}(M)$ denote the group of all
automorphisms of $M$ equipped with the 
topology of pointwise norm convergence
on the predual of $M.$ 
It was shown in \cite{ConnesT}
that 
the normal subgroup
${\rm Inn}(M)$ of inner automorphisms
of $M$ is open in ${\rm Aut}(M).$ 
Such a result cannot be expected
for $C^*$-algebras.
Indeed,  let $A$ be a separable unital $C^*$-algebra $A,$ 
and let
${\rm Aut}(A)$ be equipped with the 
topology of pointwise convergence.
It is known that, if $A$ is does not have
a continuous trace, then $\overline{{\rm Inn}(A)}/{\rm Inn}(A)$
is uncountable (\cite{Phil}). In particular,
this is true for $A=C^*_r(\Ga)$
when $\Ga$ is a countable group which is
not of type I. This is the case for ``most"
discrete groups: by \cite{Thoma}, $\Ga$ is
of type I if and only if $\Ga$ contains
an abelian normal subgroup
of finite index. However,
imitating the proof of \cite{ConnesT},
we have the following result.
Let $\Ga$ be a discrete group with
\PT. The subgroup ${\rm Inn}(L(\Ga))\cap {\rm Aut}(C^*_r(\Ga))$
is open in ${\rm Aut}(C^*_r(\Ga)).$
Observe that 
${\rm Inn}(L(\Ga))\cap {\rm Aut}(C^*_r(\Ga))$
contains $\overline{{\rm Inn}(C^*_r(\Ga))}$ when
$\Ga$ is ICC and has \PT; see \cite[Corollary 3.2 ]{Phil}.
\end{remark}

\noindent
{\bf Address}

\noindent
 Bachir Bekka, UFR Math\'ematique, Universit\'e de  Rennes 1, 
Campus Beaulieu, F-35042  Rennes Cedex, France

\noindent
E-mail : bachir.bekka@univ-rennes1.fr

\end{document}